\documentclass[12pt]{article}

\usepackage{arxiv}
\usepackage{graphicx}

\usepackage{amssymb}

\usepackage{color}

\usepackage{txfonts}
\usepackage{caption}
\usepackage{subcaption}
\newtheorem{theorem}{Theorem}
\newtheorem{proof_theorem}{Proof of Theorem}
\newtheorem{lemma}{Lemma}
\newtheorem{proof_lemma}{Proof of Lemma}

\title{Algebraic Geometrical Analysis of Metropolis Algorithm\\ When Parameters Are Non-identifiable}

\author{
    {Kenji~Nagata}\thanks{nagata.kenji@nims.go.jp} \\
    Center for Basic Research on Materials\\
	National Institute for Materials Science\\
    Tsukuba, Ibaraki 305-0044, Japan\\
	\texttt{nagata.kenji@nims.go.jp} \\
	\And
    {Yoh-ichi~Mototake}\thanks{y.mototake@r.hit-u.ac.jp} \\
    Graduate School of Social Data Science\\
	Hitotsubashi University\\
	Kunitachi, 186-8601, Tokyo, Japan\\
	\texttt{y.mototake@r.hit-u.ac.jp} \\
}

\begin{document}
\maketitle

\begin{abstract}
The Metropolis algorithm is one of the Markov chain Monte Carlo (MCMC) methods that realize sampling from the target probability distribution.
In this paper, we are concerned with the sampling from the distribution in non-identifiable cases that involve models with Fisher information matrices that may fail to be invertible.
The theoretical adjustment of the step size, which is the variance of the candidate distribution, is difficult for non-identifiable cases. 
In this study, to establish such a principle, the average acceptance rate, which is used as a guideline to optimize the step size in the MCMC method, was analytically derived in non-identifiable cases. 
The optimization principle for the step size was developed from the viewpoint of the average acceptance rate. 
In addition, we performed numerical experiments on some specific target distributions to verify the effectiveness of our theoretical results.
\end{abstract}


\section{Introduction}
The Markov chain Monte Carlo (MCMC) method~\cite{metropolis1953equation, hastings1970monte}  
is a well-known method for sampling the desired probability distribution 
and is widely used in various fields such as physics, biology, and computer science. 
In this paper, we are concerned with the sampling from the distribution in non-identifiable cases that involve models with Fisher information matrices that may fail to be invertible. 
For example, upon breaking the parameter identifiability, the Fisher information matrix becomes singular for hierarchical models, such as Gaussian mixture models~\cite{aoyagi2012learning, yamazaki2003singularities}, hidden Markov models~\cite{yamazaki2005algebraic}, and neural networks~\cite{aoyagi2012learning}.\par
Among MCMC methods, the most basic algorithm is the Metropolis algorithm. 
The procedure of the Metropolis algorithm is to set an initial value, 
generate random numbers from the current point according to the candidate distribution, 
and accept the candidate according to a defined probability.
This method has the advantage that it can be applied to any probability distribution,
regardless of the parameter identifiability.
There are various tuning parameters to improve the efficiency of the Metropolis algorithm.
In particular, the width of the candidate distribution, 
which we will call the step size hereafter, 
has a significant impact on the efficiency of the Metropolis algorithm.
However, the optimal step size depends 
on the shape, support, and width of the target distribution, 
which is a bottleneck in the application of the Metropolis algorithm.

Various studies have been conducted to theoretically realize such an optimal design.
For example, studies have been conducted to develop the optimal design 
using the ratio of candidates generated by the candidate distribution 
that is accepted (acceptance rate) as an index.
The conclusion made from these studies is that the optimal acceptance rate is 0.234 
and that this acceptance rate can be achieved by setting 
the variance of the candidate distribution 
to $O(1/d)$ with respect to the dimension $d$ 
of the sampling space~\cite{garthwaite2016adaptive}.
In addition, several automatic parameter optimization methods 
have been developed on the basis of this theory~\cite{garthwaite2016adaptive}.
On the other hand, these are theoretical analyses 
that impose a regularity condition, in the case when the parameters are identifiable, on the target distribution~\cite{roberts2001optimal}. 
Thus, the development of a theory of optimal design 
that removes the regularity assumption is expected 
to markedly expand the applicability of the Monte Carlo method.

The purpose of this study is to theoretically derive 
the relationship between step size and acceptance rate 
in the Metropolis algorithm from a general target distribution 
including non-identifiable cases.
To realize this purpose, 
a theoretical approach to models is necessary, 
in which the regularity condition does not hold.
Previously, Watanabe proposed a framework for analytically 
deriving the marginalization of a target distribution 
for which the regularity condition does not hold 
using an algebraic geometry technique 
called blowing up~\cite{watanabe2001algebraic_comp, watanabe2001algebraic}.
Furthermore, with this framework, 
an analysis of the exchange rate of the replica exchange Monte Carlo method~\cite{hukushima1996exchange}, 
a type of extended ensemble method, has been realized~\cite{nagata2008asymptotic}.
By applying this framework, we theoretically derive the relationship 
between the step size of Metropolis sampling and the acceptance rate 
from a general target distribution that includes non-identifiable cases.
Then, we performed numerical experiments on some specific target distributions to verify the effectiveness of our theoretical results.

\section{Background}

In this paper, the variable to be sampled is set to be $w \in {\mathbb R^d}$. 
Also, sampling from an exponential type probability density function of the following form is assumed,
\begin{eqnarray}
p(w) & = & \frac{1}{Z} \exp \left( - n f(w)\right) \varphi(w), \label{eq_z} \\
Z & = & \int dw \exp(-nf(w))\varphi(w),
\end{eqnarray}
where $n > 0$ and $n \in \mathbb{R}$. 
For example, in Bayesian inference, $n$ corresponds to the number of data. 
In statistical physics, $n$ corresponds to the inverse temperature. 
The function $\varphi(w)$ is defined as a probability distribution and represents a measure. 
Furthermore, $\varphi(w)$ corresponds to the prior distribution in the context of Bayesian inference. 

\subsection*{Metropolis algorithm}
The Metropolis algorithm is one of the most basic importance sampling methods based on the Markov chains algorithm. 
The sampling series can be used to approximate the target distribution and to perform integral calculations such as expectation value estimation. 
The procedure of the Metropolis algorithm is described as follows. 
\begin{enumerate}
    \item Set the parameter $w$ to the initial value $w_{(0)}$.
    \item Make a candidate $w'$ from the current state $w_{(t)}$ on the basis of the probability $q(w' | w)$.
    \item Determine the next state $w_{(t+1)}$ according to the following conditions:
    \begin{eqnarray}
        w_{(t+1)} &=&   
        \left\{
            \begin{array}{l}
                w'\:\:\:\:\:\:\:{\rm with\:\:probability\:\:} u(w', w_{(t)})  \\
                w_{(t)}\:\:\:\:\:{\rm otherwise}
            \end{array}
        \right.,\\
        u(w', w_{(t)}) &=& \min\left(1, r(w', w_{(t)})\right) \label{rate},\\
        r(w', w_{(t)}) &=& \frac{p(w_{(t)})}{p(w')},
    \end{eqnarray}

    \item Change the sampling step $t \rightarrow t+1$ and return to Step 2.
\end{enumerate}
By repeating the above procedure until $t=T$, we obtain $T$ samples. 
From the sampling series obtained by this procedure, an approximation of the target distribution or integral calculations is realized. 
The probability $q(w' | w)$ of generating a candidate $w'$ is called the candidate distribution. 
When the detailed balancing condition is satisfied, the sample series generated by the Metropolis-Hastings algorithm can be regarded as a sample from the probability distribution $p(w)$. 
The detailed balance condition is realized by assuming that the candidate distribution has the property $q(w' | w) = q(w | w')$. 

\subsection*{Parameter design of Metropolis algorithm}
To run the Metropolis algorithm appropriately, it is necessary to optimize the width of the candidate distribution. 
We will refer to this width as ``step size''. 
To achieve efficient sampling from a probability distribution, it is desirable for the samples in a sampling series to be uncorrelated with each other. 
To achieve this, the step size should be set larger.  
However, if the step size is set to large, the probability $r(w', w_{(t)})$ in Eq.~(\ref{rate}) will be significantly smaller and the parameter $w_{(t)}$ becomes hardly updated. 
If the step size is reduced, the probability $r(w', w_{(t)})$ tends to become larger, but it is difficult to cover the space that should be sampled because the sampling series is limited to only local regions. 
Therefore, it is important to optimize the step size to achieve efficient sampling. 
In the optimization, the criterion that the acceptance ratio reaches an appropriate value is used. 
\par

\section{Main Analysis: derivation of average acceptance rate}
In this section, we present a theoretical analysis of the average acceptance rate in the Metropolis algorithm and discuss its relationship with step size and target distribution.

\subsection*{Problem Setting}
Suppose that $w \in {\mathbb R^d}$, we set the target distribution $p(w)$ as follows:
\begin{eqnarray}
p(w) & = & \frac{1}{Z} \exp \left( - n f(w)\right) \varphi(w), \\
&&\forall w: f(w)  \geq  0, \: f(0) = 0, \label{target_pd1}
\end{eqnarray}
where $Z$ is defined by Eq.~(\ref{eq_z}).
Note that Eq.~(\ref{target_pd1}) does not indicate that $f(w) = 0$ is only satisfied at $w=0$. 
Suppose that the probability distribution $\varphi(w)$ has  a bounded closed set containing $w=0$ as its support.
The target distribution $p(w)$ does not necessarily converge to a Gaussian distribution with $n\rightarrow \infty$.
This is due to the condition that the Hessian of $f(w)$ is not necessarily a positive definite matrix.
In other words, this analysis also considers the case where singularities are included in the algebraic set defined by $f(w)=0$.
Here, we assume that the singularity of the algebraic set $f(w)=0$ is at the origin $w=0$.
We consider performing the Metropolis algorithm on this target distribution $p(w)$.
The candidate distribution $q(w' | w)$ is set as a one-dimensional Gaussian distribution as follows:
\begin{eqnarray}
q(w' | w ) & = & \mathcal{N}(w'_i;w_i,\sigma^2) \prod_{j\neq i}\delta(w'_j - w_j). \label{candidate}
\end{eqnarray}
Then, the acceptance rate $u(w,w')$ is given by
\begin{eqnarray}
u(w,w') & = & \min(1,r(w,w')), \label{accept}\\
r(w,w') & = & \frac{p(w')}{p(w)}.
\end{eqnarray}
By using this acceptance rate, we define the average acceptance rate as the following expected calculation:
\begin{eqnarray}
U & = & \int dw p(w) \int dw' \: u(w,w') q(w' | w). \label{ave.accept}
\end{eqnarray}
The main task of this study is to analytically determine the acceptance rate $U$.
In preparation for this analysis, we define the following function $Z_i$.
\begin{eqnarray}
Z_i & = & \int dw |w_i| \exp(-nf(w))\varphi(w) 
\end{eqnarray}
We also define two types of zeta functions, $\zeta(z)$ and $\zeta_i(z)$, 
for the function $f(w)$ and the probability distribution $\varphi(w)$ as follows:
\begin{eqnarray}
\zeta(z) & = & \int dw f(w)^z \varphi(w), \\
\zeta_i(z) & = & \int dw \:|w_i| \: f(w)^z \varphi(w), 
\end{eqnarray}
where $z$ is a one-dimensional complex number.
The zeta functions $\zeta(z)$ and $\zeta_i(z)$ are regular functions in the domain $Re(z)>0$ 
and can be analytically connected to rational-type functions in the entire complex domain.
All its poles are rational numbers on the negative real axis.
We define the rational number $-\lambda$ as the pole closest to the origin of the zeta function $\zeta(z)$, 
and define the natural number $m$ as its order.
If Hessian is positive definite for any $w$, then $\lambda=\frac{d}{2}$ and $m=1$ are satisfied.
Otherwise, $\lambda$ and $m$ can be analytically calculated 
by using the resolution of singularity in algebraic geometry~\cite{atiyah1970resolution}.
Similarly, we define the poles and its order 
corresponding to the zeta function $\zeta_i(z)$ as $-\lambda_i$ and $m_i$, respectively.
\subsection*{Main Theorem}
\begin{theorem}
\label{theorem1}
The average acceptance rate $U$ has the following asymptotic form
with a sufficiently large step size $\sigma$ in $n\rightarrow\infty$,
\begin{eqnarray}
U & \sim & \frac{(\log n)^{m_i-m}}{n^{\lambda_i-\lambda}} \frac{c 4\sqrt{2}}{\sigma}, \label{U_general}
\end{eqnarray}
where $c$ is a constant determined by the target distribution $p(w)$ and the candidate distribution $q(w'|w)$.
\end{theorem}

To derive theorem \ref{theorem1}, we first show the following Lemma 1,
\begin{lemma}
\label{lemma_1}
The average acceptance rate $U$ is expressed using the functions $Z$ and $Z_i$ as follows:
\begin{eqnarray}
U & \sim & \frac{4\sqrt{2}}{\sigma} \frac{Z_i}{Z}.
\end{eqnarray}
\end{lemma}

\begin{proof_lemma}
\label{proof_lemma_1}
From Eqs.~(\ref{candidate}), (\ref{accept}), (\ref{ave.accept}) and the symmetry property of the candidate distribution $q(w'|w)$, $q(w'|w)=q(w|w')$, the average acceptance rate $U$ is expressed as
\begin{eqnarray}
U  & = & \int dw p(w) \left\{\int_{p(w') > p(w)} dw'   q(w' | w) + \int_{p(w') < p(w)} dw'  \: r(w,w') q(w' | w) \right\} \nonumber \\
 & = & \int dw \left\{\int_{p(w') > p(w)} dw'  p(w) q(w' | w) + \int_{p(w') < p(w)} dw'  \: p(w') q(w| w') \right\} \nonumber \\
 & = & 2 \int dw p(w) \int_{p(w') > p(w)} dw'  q(w' | w) \label{p_p}\\
 & \sim & \frac{2}{Z} \int dw \exp(-nf(w))  \varphi(w) \int_{f(w') < f(w)} dw' \:  q(w' | w).
 \label{transition}
\end{eqnarray}
Here, the set $\{ w' ; p(w') > p(w)\}$, defined as the integral interval of $w'$ in Eq. (\ref{p_p}),
is approximated using $\{ w' ; f(w') < f(w)\}$ on the basis of the fact that $n$ is sufficiently large.

By substituting Eq. (\ref{transition}), we get
\begin{eqnarray}
U & = & \frac{2}{Z} \int dw \exp(-nf(w))  \varphi(w) \\
&& \times \int_{f(w') < f(w)} dw' \mathcal{N}(w'_i;w_i,\sigma^2) \prod_{j\neq i}\delta(w'_j-w_j)\\
& \sim & \frac{2}{Z} \int dw \exp(-nf(w))  \varphi(w)\int^{|w_i|}_{-|w_i|} dw'_i  \mathcal{N}(w'_i;w_i,\sigma^2)\\
& = & \frac{\tilde{Z}}{Z}. \label{Z_Z}
\end{eqnarray}
Here, on the basis of the fact that the singularity of the set $f(w)=0$ is at the origin, 
we approximate the integral range $\{ w'_i ; f(w_1,\cdots,w'_i,\cdots,w_d) < f(w_1,\cdots,w_i,\cdots,w_d)\}$ 
as $\{w'_i ; -|w_i| < w'_i < |w_i |\}$.
The function $\tilde{Z}$ can be expressed as follows by integrating over $w_i$.
\begin{eqnarray}
\tilde{Z} & = & 2 \int dw \exp(-nf(w)) \varphi(w) \int^{|w_i|}_{-|w_i|} dw'_i  \mathcal{N}(w'_i;w_i,\sigma^2) \\
 &=& 2\int dw \exp(-nf(w)) \varphi(w) \sqrt{\pi}\: {\rm erf}\left(\frac{\sqrt{2}|w_i|}{\sigma}\right)
\end{eqnarray}
A Taylor expansion of the ${\rm erf}$ function is expressed as follows:
\begin{eqnarray}
{\rm erf}\left(\frac{\sqrt{2}|w_i|}{\sigma}\right) \sim \frac{2}{\sqrt{\pi}}\frac{\sqrt{2}|w_i|}{\sigma}.
\end{eqnarray}
By using this, we obtain the function $\tilde{Z}$ is as follows:
\begin{eqnarray}
\tilde{Z} &=& \frac{4\sqrt{2}}{\sigma}\int dw \: |w_i|\:\exp(-nf(w)) \varphi(w) \\
 & = & \frac{4\sqrt{2}}{\sigma} Z_i.
\end{eqnarray}
By substituting this equation into the equation (\ref{Z_Z}), we get the following equation:
\begin{eqnarray}
U &=&  \frac{4\sqrt{2}}{\sigma} \frac{Z_i}{Z},
\end{eqnarray}
which completes the proof of Lemma 1.
\begin{flushright}
$\Box$
\end{flushright}
\end{proof_lemma}

From Lemma \ref{proof_lemma_1}, the asymptotic form of the average acceptance rate $U$ 
can be obtained by deriving the functions $Z$ and $Z_i$.
Then, we show Lemma 2.

\begin{lemma}
\label{lemma_2}
The functions $Z$ and $Z_i$ respectively have the following asymptotic forms under $n \rightarrow \infty$:
\begin{eqnarray}
Z & = & \frac{(\log n)^{m-1}}{n^\lambda} g(-\lambda)\Gamma(\lambda), \\
Z_i & = & \frac{(\log n)^{m_i-1}}{n^{\lambda_i}} g_i(-\lambda_i) \Gamma(\lambda_i).
\end{eqnarray}
Here, $g(-\lambda)$ and $g_i(-\lambda_i)$ are regular functions and $\Gamma(\lambda)$ is the gamma function.
\end{lemma}

\begin{proof_lemma}
\label{proof_lemma_2}
Refer to \cite{watanabe2001algebraic} for the derivation of the asymptotic form of the function $Z$.
In this proof, we derive an asymptotic form for $Z_i$.
The function $Z_i$ is expressed using the delta function $\delta(\cdot)$ as follows:
\begin{eqnarray}
Z_i &=& \int dw \int_{0}^{\infty} dt \delta(t - f(w)) \:|w_i|\:\exp(-nt) \varphi(w)\\
&=& \int_{0}^{\infty} dt \:\exp(-nt) \int dw  \delta(t - f(w))\:|w_i|\:\varphi(w)\\
& = & \int_{0}^{\infty} dt \exp(-nt) V_i(t),
\label{EQ_laplace}
\end{eqnarray}
where the function $V_i(t)$ is defined as
\begin{eqnarray}
V_i(t) = \int dw \delta(t - f(w)) \:|w_i|\:\varphi(w).
\end{eqnarray}
Equation~(\ref{EQ_laplace}) shows that the Laplace transform for $V_i(t)$ is equal to $Z_i$.
By applying the Mellin transform to $V_{i}(t)$, we obtain
\begin{eqnarray}
\zeta_i(z) &=& \int_{0}^{\infty} dt \: t^z V_i(t)\\
&=& \int_{0}^{\infty} dt \: t^z \int dw \delta(t - f(w)) \:|w_i|\:\varphi(w) \\
&=& \int dw \:|w_i|\:f(w)^z \varphi(w). \label{zeta'}
\end{eqnarray}
This function has poles only on the negative real axis, 
and the pole closest to the origin is known to have a major contribution to this function~\cite{watanabe2001algebraic}.
Now, if the position of this pole is $-\lambda_i$ and its order is $m_i$, then $\zeta_i(z)$ is expressed by
\begin{eqnarray}
\zeta_i(z) \sim \frac{g_i(z)}{(z + \lambda_i)^{m_i}},
\end{eqnarray}
where $g_i(z)$ is defined to be a regular function with no poles.
By inverse Mellin transformation of this function, we can derive the function $V_{i}(t)$ as follows:
\begin{eqnarray}
V_i(t) & = & \frac{1}{2 \pi {\rm i}} \int_{c - {\rm i}\infty}^{c + {\rm i}\infty} \zeta_i(z) t^{-z-1} dz \\
 & = & \frac{1}{2 \pi {\rm i}} \int_{c - {\rm i}\infty}^{c + {\rm i}\infty} \frac{g_i(z) t^{-z-1}}{(z + \lambda_i)^{m_i}} dz \\
 & = & \left.\frac{\partial^{m_i-1}}{\partial z^{m_i-1}} \left\{g_i(z) t^{-z-1} \right\} \right|_{z=-\lambda_i} \\
 & = & g_i(-\lambda_i) t^{\lambda_i-1} (- \log t)^{m_i-1} \left( 1 + O\left( \frac{1}{\log t}\right)\right) \\
 & \sim & g_i(-\lambda_i) t^{\lambda_i-1} (- \log t)^{m_i-1}.
\end{eqnarray}
Therefore, $Z_i$ can be expressed by
\begin{eqnarray}
Z_i & = & \int_{0}^{\infty} dt \exp(-nt) V_{i}(t) \\
& = & \int_{0}^{\infty} dt \exp(-nt) g_i(-\lambda_i) t^{\lambda_i-1} (- \log t)^{m_i-1} \\
 & = & g_i(-\lambda_i) \int_{0}^{\infty} \frac{dt'}{n} \exp(-t') \left( \frac{t'}{n}\right)^{\lambda_i-1} (\log n - \log t')^{m_i-1} \\
 & = & \frac{(\log n)^{m_i-1}}{n^{\lambda_i}} g_i(-\lambda_i) B(\lambda_i,m_i). 
\end{eqnarray}
Here, $B(\lambda_i,m_i)$ can be calculated as follows:
\begin{eqnarray}
B(\lambda_i,m_i) & = & \int_{0}^{\infty} dt \exp(-t) t^{\lambda_i-1} \left(1 - \frac{\log t'}{\log n}\right)^{m_i-1} \sim \Gamma (\lambda_i).
\end{eqnarray}
Consequently, $Z_i$ can be expressed by
\begin{eqnarray}
Z_i & = & \frac{(\log n)^{m_i-1}}{n^{\lambda_i}} g_i(-\lambda_i)\Gamma (\lambda_i), \label{result_Z_i}
\end{eqnarray}
which completes the proof of Lemma 2.
\begin{flushright}
$\Box$
\end{flushright}
\end{proof_lemma}

Finally, we prove Theorem 1 by using Lemma 1 and Lemma 2.

\begin{proof_theorem}
By substituting $Z$ and $Z_i$ derived in Lemma 2 into the result of Lemma 1, we obtain Theorem 1.
\begin{eqnarray}
U & \sim & \frac{4\sqrt{2}}{\sigma} \frac{Z_i}{Z} \\
 & = & \frac{4\sqrt{2}}{\sigma} 
    \frac{(\log n)^{m_i-m}}{n^{\lambda_i - \lambda}}
    \frac{g_i(-\lambda_i)\Gamma (\lambda_i)}{g(-\lambda)\Gamma (\lambda)}\\
 & = & \frac{(\log n)^{m_i-m}}{n^{\lambda_i - \lambda}}
    \frac{A(\lambda,\lambda_i) 4\sqrt{2}}{\sigma}
\end{eqnarray}
Here, the function $A(\lambda,\lambda_i)$ is defined by 
\begin{eqnarray}
A(\lambda,\lambda_i) & = & \frac{g_i(-\lambda_i)\Gamma (\lambda_i)}{g(-\lambda)\Gamma (\lambda)}.
\end{eqnarray}
\begin{flushright}
$\Box$
\end{flushright}
\end{proof_theorem}

The expression for the average acceptance rate (Eq.~(\ref{U_general})) is obtained as an asymptotic form when the number of data $n$ is sufficiently large. 
The property of the average acceptance rate can be categorized by a combination of poles $\lambda$ and $\lambda_i$ and orders $m$ and $m_i$. 
The combination pattern of poles and orders, which characterizes the behavior of the average acceptance rate, can be represented by the following three cases,
\begin{description}
    \item[Case 1:] $\lambda < \lambda_i$,
    \item[Case 2:] $\lambda = \lambda_i$,$m > m_i$,
    \item[Case 3:] $\lambda = \lambda_i$,$m = m_i$,
\end{description}
where, we do not consider the case of $\lambda > \lambda_i$ and the case of $\lambda = \lambda_i$, $m < m_i$. 
The reason why we limit to only such cases is that the integrand function of $\zeta_i(z)$ is defined as the integrand function of $\zeta(z)$ multiplied by $|w_i|$. 
From this, there is no need to consider the case where $\lambda_i$ is smaller than $\lambda$, or the case where $\lambda$ = $\lambda_i$ and the order decreases. 

\subsection*{Limitation of Theorem 1}
We discuss two points about the limitation of Theorem 1. 
First, we discuss the limitation of the theorem on step size. 
In the proof of the theorem, we approximated equation (\ref{Z_Z}) by assuming that the candidate distribution contains a singular point in the neighborhood of the origin. 
Therefore, in regions with small step sizes, the acceptance rate derived from the theorem and the actual acceptance rate will deviate. 
The range of step sizes for which the theorem holds depends on the constants $n$, $\lambda$, and $m$. 
The exact scope of the theorem is a subject for future work.\par
Second, we discuss the limitation of the theorem on constant $n$. 
Owing to the property of asymptotic expansion, $\frac{(\log n)^{m_i-m}}{n^{\lambda_i-\lambda}}$ becomes the main term of the acceptance rate when $n$ is sufficiently  large~\cite{watanabe2001algebraic_comp, watanabe2001algebraic}. 
On the other hand, the accuracy of terms of the lower order than the main term is not guaranteed. 
Therefore, even if it could be derived theoretically, the constant term $c$ theoretically obtained would not match with the simulation result. 
Therefore, in the verification of the theorem by simulation, which will be explained later, only the main term will be compared.

\renewcommand{\arraystretch}{1.5}
\begin{table}[tbh]
\centering
\caption{List of parameters $\lambda$ and $m$ for two cases of numerical simulation.}
\begin{tabular}{|c|c|c|c|c|c|c|}
\hline 
case & $\lambda$ & $m$ & $\lambda_1$ & $m_1$ & $\lambda_2$ & $m_2$ \\
\hline 
$f(w) = w^2_1w^2_2$ & $1/2$ & 2 & $1/2$ & 1 & $1/2$ & 1  \\
\hline 
$f(w) = w^2_1w^4_2$ & $1/4$ & 1 & $1/4$ & 1 & $1/2$ & 2  \\
\hline 
\end{tabular}
\label{tbl1}
\end{table}
\renewcommand{\arraystretch}{1.0}

\section{Numerical Experiments}
In this section, we describe a numerical experiment to verify Theorem 1. 
We conducted numerical experiments for the following two types of function $f(w)$; 
\begin{itemize}
    \item[(a)] $f(w) = w^2_1w^2_2$, \label{case1}
    \item[(b)] $f(w) = w^2_1w^4_2$. \label{case2}
\end{itemize}
When $f(w)$ is set as shown above, $\lambda$, $\lambda_1$, $\lambda_2$, $m$, $m_1$, and $m_2$ are given as in Table \ref{tbl1}.
The probability distribution $\varphi(w)$ is assumed to be a two-dimensional standard normal distribution. 

As can be seen from Table \ref{tbl1}, to verify Theorem 1 in all $\mathbf{ cases\:1,\:2 and\:3}$, it is sufficient to compare Theorem 1 with the average acceptance ratio of the parameters $w_1$ and $w_2$ of these two functions. 
Specifically, the sampling of $w_1$ and $w_2$ in function (a) corresponds to $\mathbf{ case\:2}$, the sampling of $w_1$ in function (b) corresponds to $\mathbf{ case\:3}$, and the sampling of $w_2$ corresponds to $\mathbf{ case\:1}$ in function (b). 
On the basis of Theorem 1, the average acceptance rates corresponding to functions (a) and (b) are derived (see Appendices A and B). 
Theorem 1 was verified by comparing the average acceptance rate derived from Theorem 1 with the average acceptance rate numerically obtained by the Metropolis algorithm. 
For accurate verification, it is desirable that the sampling by the Metropolis-Hastings algorithm achieves sampling from the stationary distribution (Eq.~\ref{eq_z}) as closely as possible. 
For this reason, in this verification, the replica-exchange Monte Carlo method, which converges quickly to a stationary distribution, was adopted as the sampling method. 
\begin{figure}[tb]
  \begin{tabular}{cc}
  \includegraphics[clip,width=6.4cm]{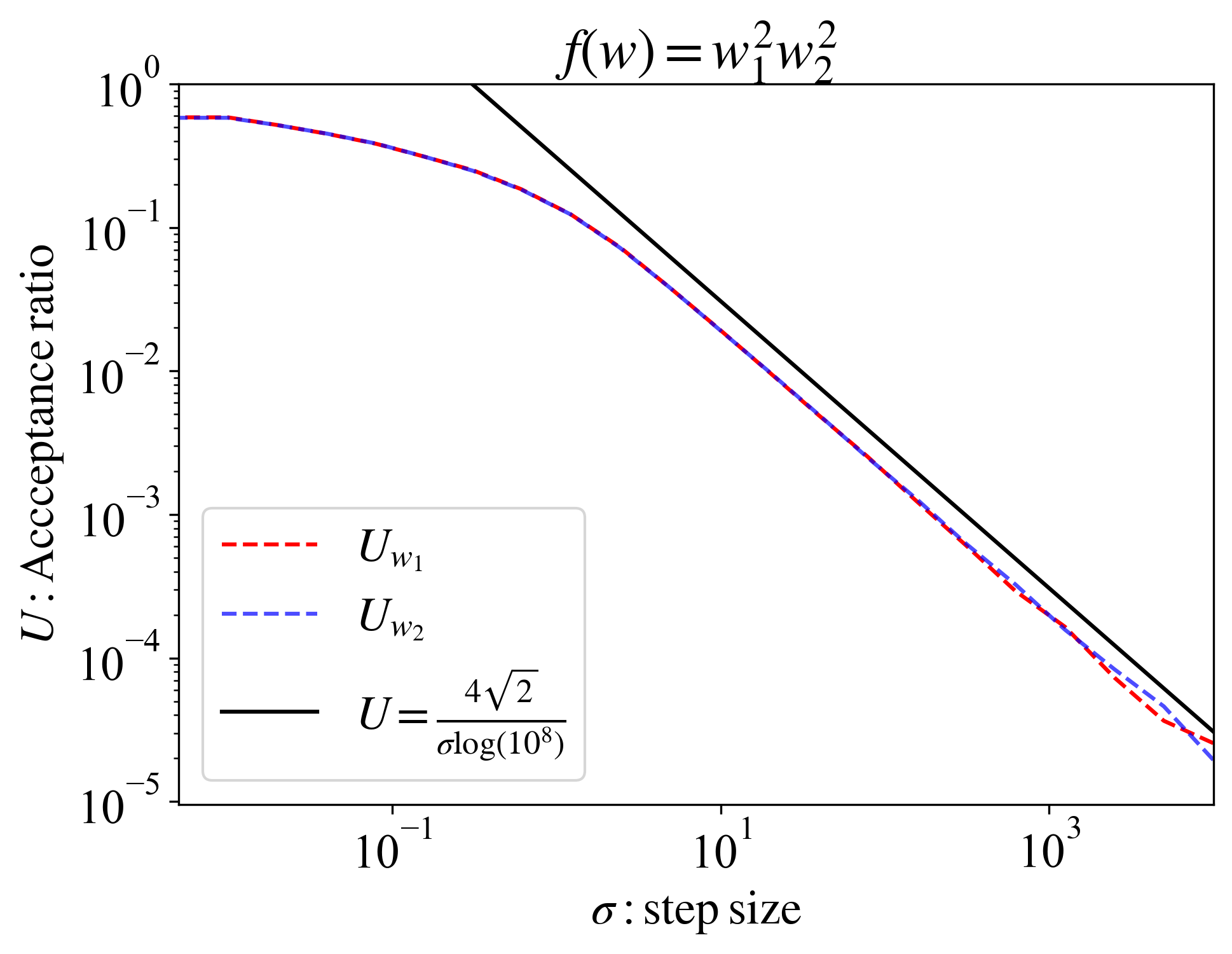} & 
  \includegraphics[clip,width=6.4cm]{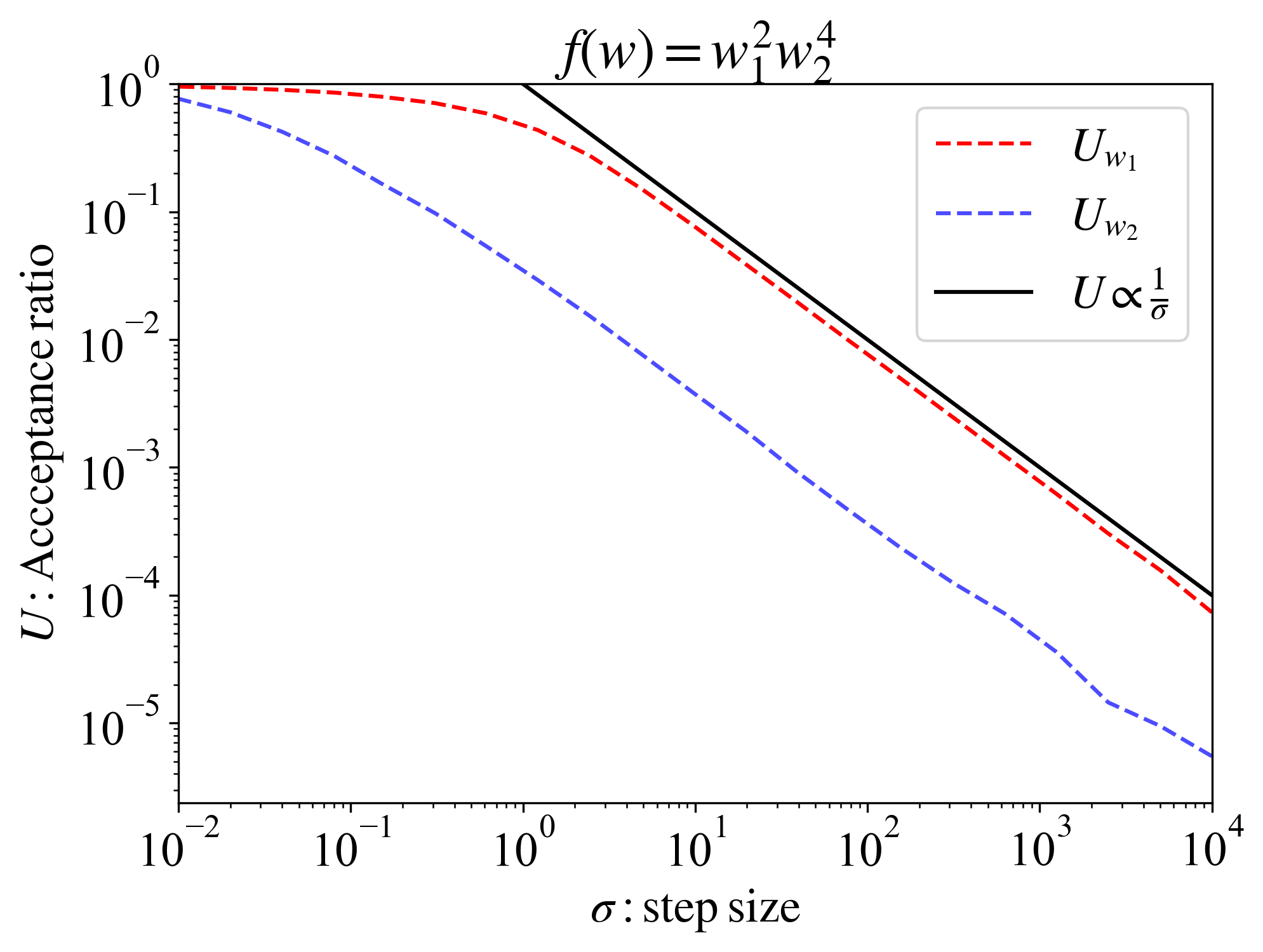}\\
  \end{tabular}
  \caption{Behavior of the average acceptance rate with respect to step size $\sigma$. Left: Relationship for $f(w)=w_1^2w_2^2$. Right: Relationship for $f(w)=w_1^2w_2^4$.}
   \label{fig_1}
\end{figure}
As can be seen from Theorem 1, the average acceptance rate is dominated by the step size and the constant $n$. 
First, we examined the relationship between the step size and the acceptance rate. 
The results of the simulation with $n=10^8$ are shown as the red and blue dotted lines in Fig.~\ref{fig_1}. 
The relationship between the step size and the average acceptance rate derived from Theorem 1, $U\propto\frac{1}{\sigma}$, is shown by the black line. 
To easily compare the theoretical and experimental behaviors of the average acceptance rate as a function of step size, the theoretical curve was corrected to match the simulated value at the right endpoint of the visible region of the graph ($\sigma=10^4$). 
The results of this comparison confirm that in the region where the step size $\sigma$ is sufficiently large, the behavior of the average acceptance rate in theory and simulation is consistent. 

\begin{figure}[tb]
  \begin{tabular}{cc}
  \includegraphics[clip,width=6.4cm]{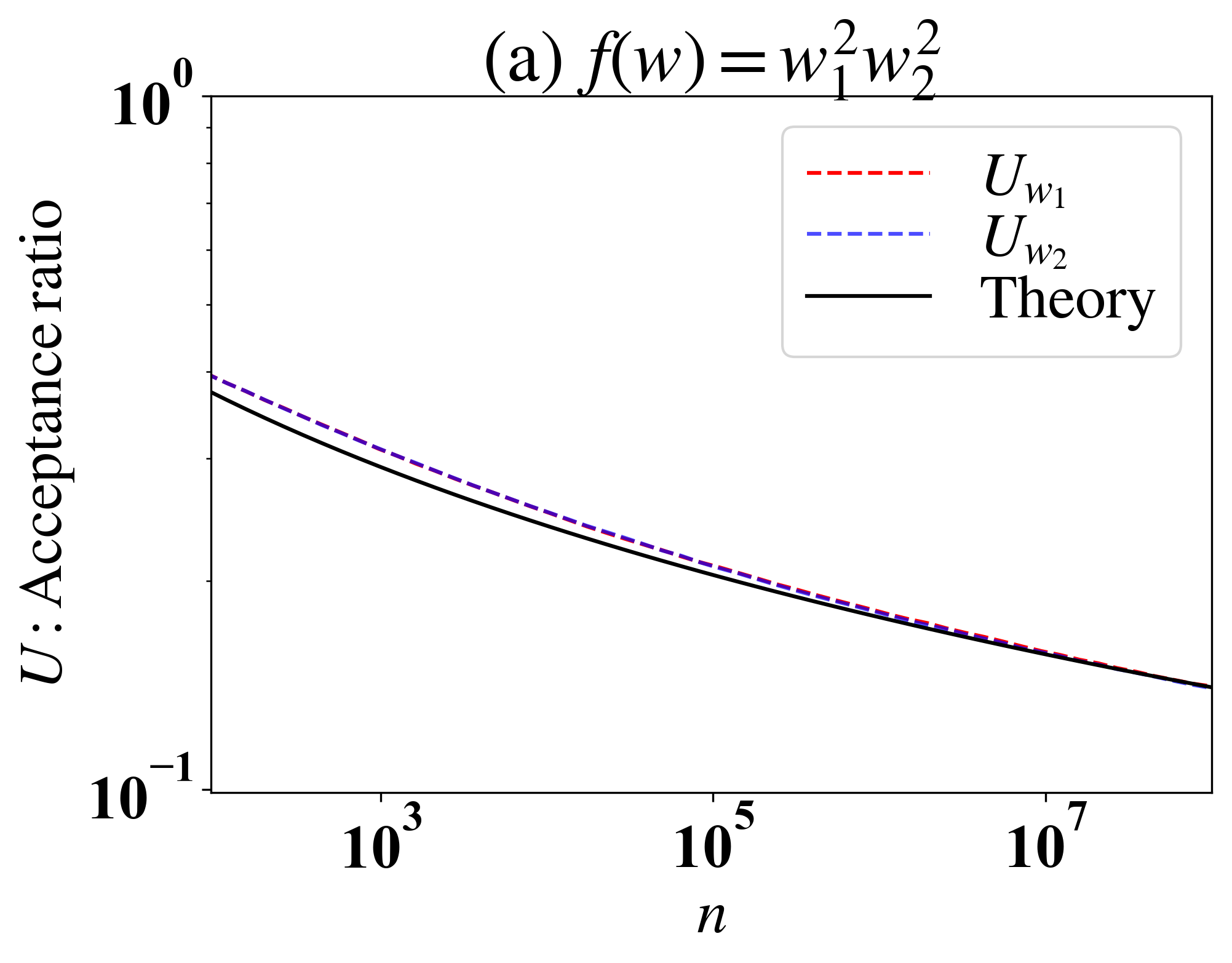} & 
  \includegraphics[clip,width=6.4cm]{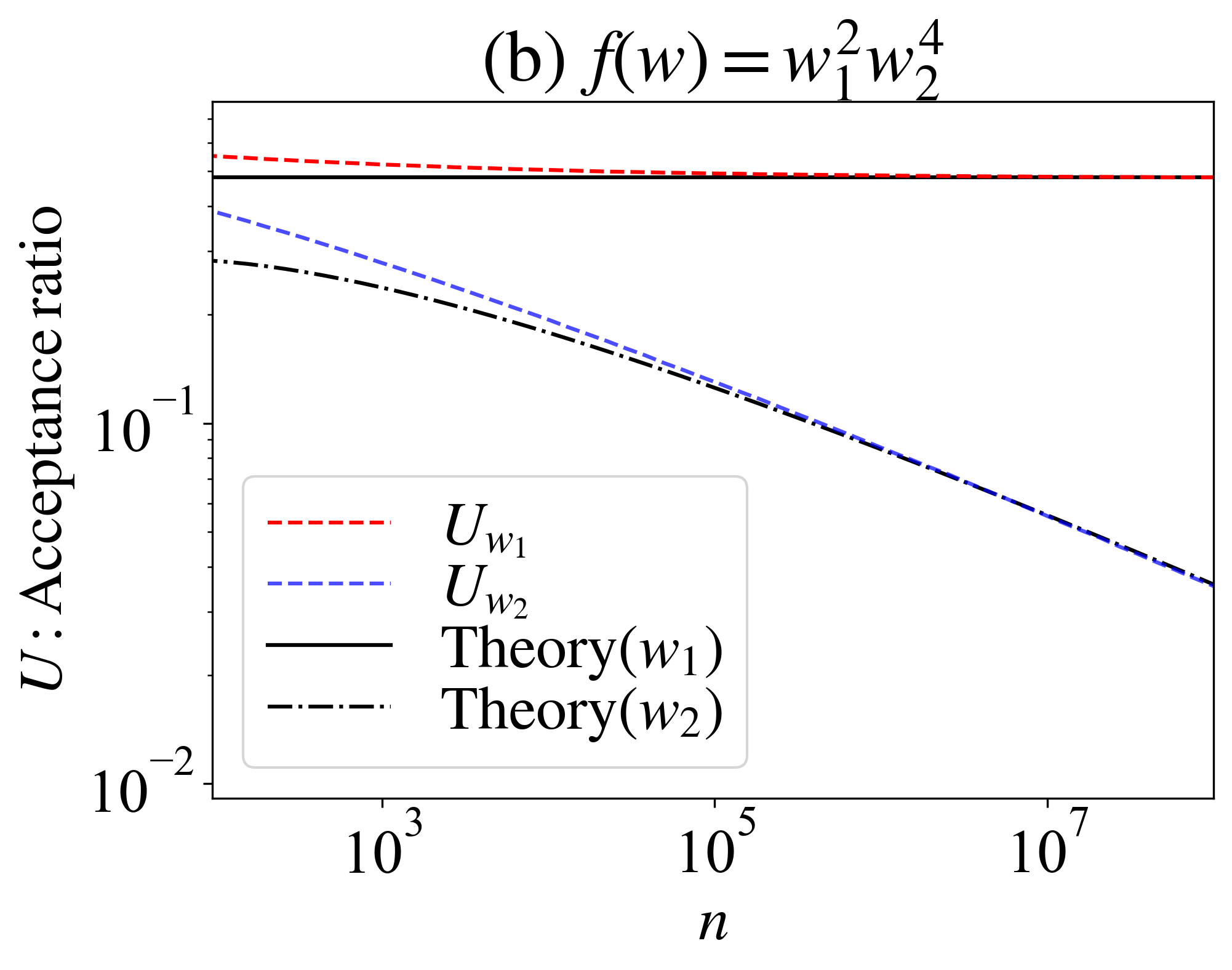}\\
  \end{tabular}
  \caption{Behavior of the average acceptance rate with respect to constant $n$. Left: Relationship for $f(w)=w_1^2w_2^2$. Right: Relationship for $f(w)=w_1^2w_2^4$.}
   \label{fig_2}
\end{figure}

Next, we verify the theorem when the step size is fixed and the constant $n$ is varied. 
The results of the simulation with $n=10^8$ are shown as the red and blue dotted lines in Fig.~\ref{fig_2}. 
The relationship between the constant $n$ and the average acceptance rate derived from Theorem 1 is shown by the black line in Fig.~\ref{fig_2}.  
As in Fig.~\ref{fig_1}, to easily compare the theoretical and experimental behaviors of the average acceptance rate as a function of constant $n$, the theoretical curve was corrected to match the simulated value at the right endpoint of the visible region of the graph. 
The results of this comparison confirm that in the region where the constant $n$ is sufficiently large, the behavior of the average acceptance rate in theory and simulation is consistent. 

\subsection*{Conditions to keep the acceptance rate constant}
Some sampling methods, such as simulated annealing~\cite{kirkpatrick1983optimization} or the replica-exchange Monte Carlo method~\cite{hukushima1996exchange}, achieve efficient sampling by combining sampling results at different constants $n$. 
In such a method, it is important to set the optimal step size at different $n$ values. 
In this section, on the basis of Theorem 1, we derive a way to change the step size so that the average acceptance rate remains constant even if $n$ changes. 
Then, by setting the step size on the basis of the derived equation, we verify whether the average acceptance rate remains constant in numerical experiments. 

The step size at which the average acceptance rate does not change for different $n$ values is derived for $\mathbf{cases 1}$, $\mathbf{2}$, and $\mathbf{3}$ as follows. 
\begin{description}
\item[case 1:] $\lambda < \lambda_i$:
In this case, from Eq. (\ref{U_general}), the step size $\sigma$ that makes the average acceptance rate constant can be achieved by the following:
\begin{eqnarray}
\label{sigma_case1}
\sigma \sim 4\sqrt{2}A(\lambda,\lambda_i) \frac{(\log n)^{m_i-m}}{n^{\lambda_i-\lambda}}. 
\end{eqnarray}
To keep the average acceptance rate constant, this relationship indicates that the step size $\sigma$ should be reduced as $n$ increases. 
\item[case 2:] $\lambda = \lambda_i$, $m > m_i$:
From the fact that $\lambda = \lambda_i$, the main term $n^{\lambda_i-\lambda}$ in Eq. (\ref{U_general}) is constantly $1$. 
On the other hand, since $m>m_i$, $\frac{1}{(\log n)^{m-m_i}}$ remains the main term. 
Therefore, the average acceptance rate can be kept constant by setting the step size $\sigma$ as follows:
\begin{eqnarray}
\label{sigma_case2}
\sigma  \sim  \frac{4\sqrt{2}A(\lambda, \lambda)}{(\log n)^{m-m_i}}.
\end{eqnarray}

\item[case 3:] $\lambda = \lambda_i$, $m = m_i$:
As in $\mathbf{ case\:2}$, since $\lambda = \lambda_i$, $n^{\lambda_i-\lambda}$ in Eq. (\ref{U_general}) is constantly $1$. 
Also, since $m=m_i$, the term for $n$ in Eq. (\ref{U_general}) vanishes. 
Therefore, the average acceptance rate can be kept constant by setting the step size $\sigma$ as follows: 
\begin{eqnarray}
\label{sigma_case3}
\sigma \sim  4\sqrt{2}A(\lambda, \lambda). 
\end{eqnarray}
\end{description}

\begin{figure}[tb]
  \begin{tabular}{cc}
  \includegraphics[clip,width=6.4cm]{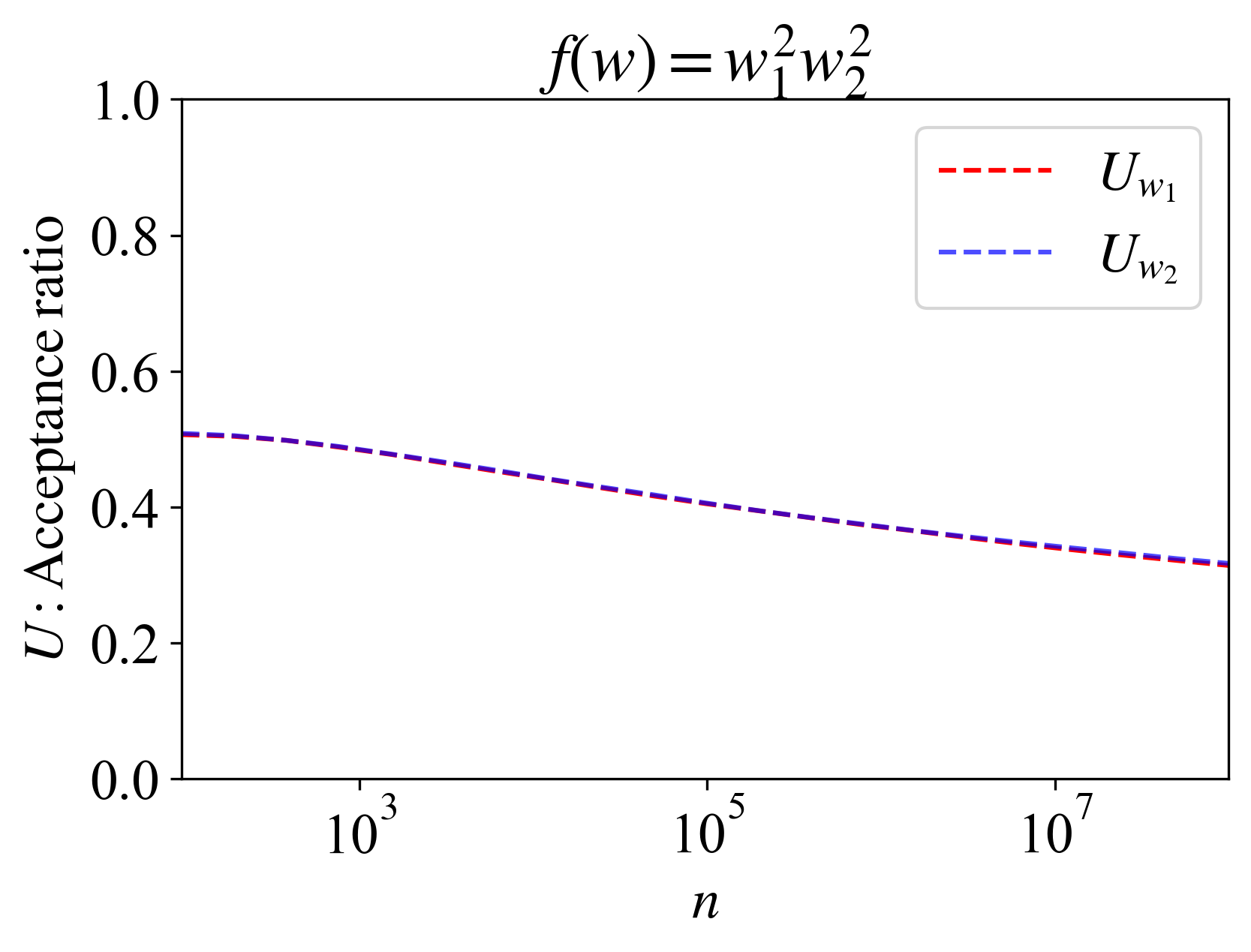} & 
  \includegraphics[clip,width=6.4cm]{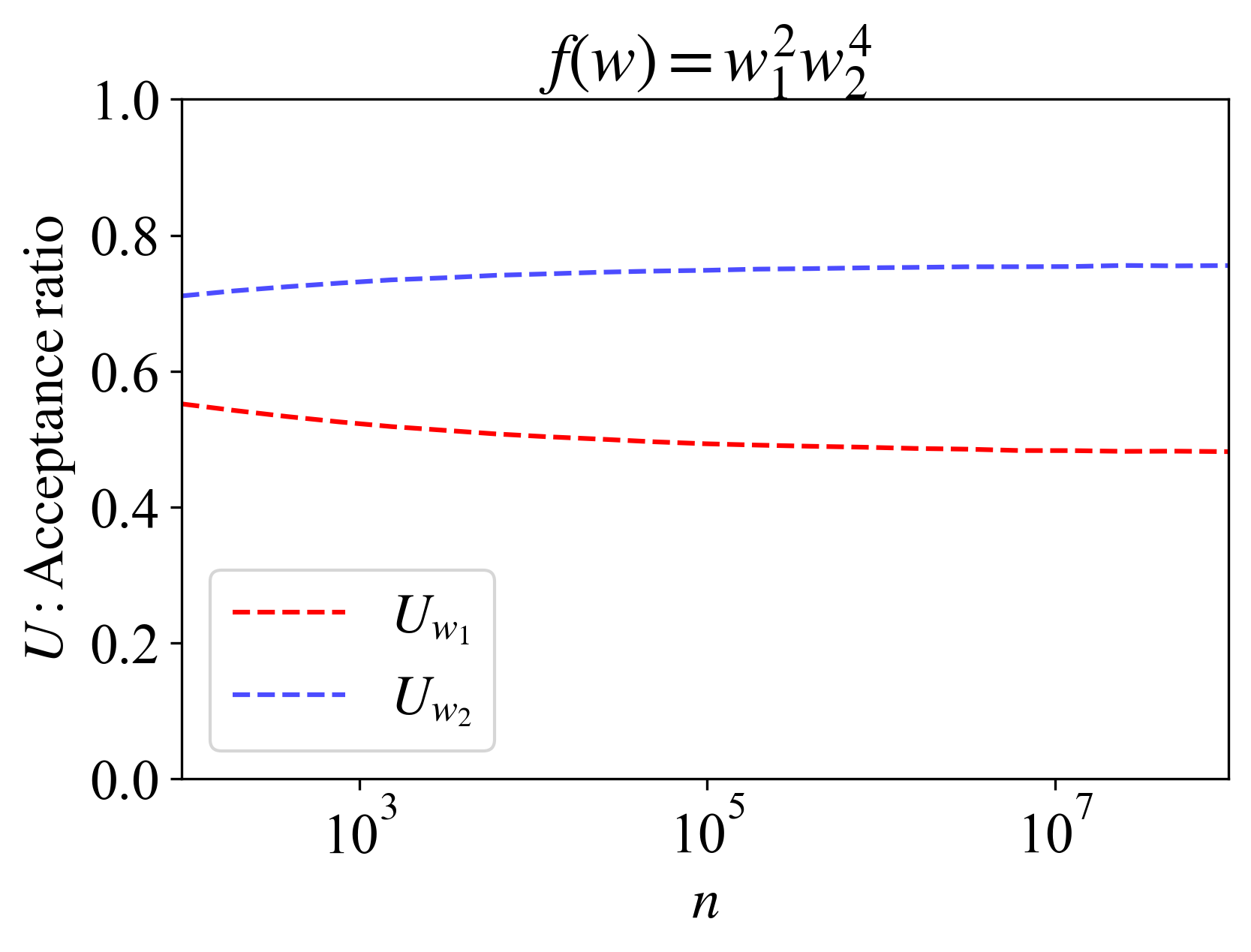}\\
  \end{tabular}
  \caption{Behavior of the average acceptance rate with respect to $n$ in the numerical experiments described in Sect. 4.1. Left: Relationship for $f(w)=w_1^2w_2^2$. Right: Relationship for $f(w)=w_1^2w_2^4$.}
   \label{fig_3}
\end{figure}

By setting the step size according to Eqs. (\ref{sigma_case1})-(\ref{sigma_case3}), we verified by numerical experiments whether the average acceptance rate remains constant independent of the constant $n$. 
The results are shown in Fig.~\ref{fig_3}. 
Note that the vertical axis in Fig.~\ref{fig_2} is a logarithmic scale, whereas the vertical axis in Fig.~\ref{fig_3} is a real scale. 
From the results, it was confirmed that the average acceptance rate is roughly constant; particularly, it is more precisely constant in the region where $n$ is large. 
This verification confirms that if $\lambda$ and $m$ are known in advance, the optimal step size for simulated annealing~\cite{kirkpatrick1983optimization} or replica-exchange Monte Carlo methods~\cite{hukushima1996exchange} can be estimated approximately on the basis of Theorem 1.

\section{Summary and Discussion}
In this study, the average acceptance rate was analytically derived to optimize the step size, 
which is important for the design of the Metropolis algorithm.
As a result, the average acceptance rate was derived as a universal expression 
using the pole $\lambda$ and its order $m$ of the zeta function defined for the target distribution.
In addition, we performed simulations of the Metropolis algorithm on a simple system, 
and compared the obtained theoretical values with the experimental values, 
which showed that the theory obtained in this study is appropriate.

The poles $\lambda$ of the zeta function are generally difficult to find in analysis.
It is possible to calculate $\lambda$ analytically 
toward an arbitrary model if the blowing-up calculation can be realized~\cite{watanabe2001algebraic_comp}.
However, the computation of blowing-up is very difficult.
For example, the Bayesian estimation framework that has been solved in practice is limited to a few stochastic models such as three-layer neural networks~\cite{aoyagi2012learning}, Gaussian mixture models~\cite{aoyagi2012learning, yamazaki2003singularities}, and hidden Markov models~\cite{yamazaki2005algebraic}.

On the other hand, it is possible to numerically derive $\lambda$ and $m$ 
from the simulation results of the Metropolis algorithm with this theorem.
By using \cite{garthwaite2016adaptive, doi:10.7566/JPSJ.90.034001}, 
an automatic optimization algorithm for step size in the Metropolis algorithm, 
we can make the average acceptance rate constant for any $n$.
By investigating the relationship between the step size obtained in this way and the constant $n$, 
we can derive $\lambda$ and $m$ numerically on the basis of this theorem.
One of the major contributions of this paper is that $\lambda$ and $m$,
which are difficult to derive analytically, can be derived numerically via the Metropolis algorithm.

To verify this theorem, we conducted simulations in cases where $\lambda$ and $m$ can be derived analytically. 
However, it is important to verify the theorem in more practical cases, and this will be a future work.

\section*{Acknowledgments}
This work was supported by KAKENHI grant numbers JP20H04648 and ISM Cooperative Research Program (2020-ISMCRP-2070).

\appendix
\section{Analytical solution for $f(w_1,w_2) = w_1^2w_2^2$}
From Lemma 2, the average acceptance rates $U_1$ and $U_2$ for $w_1$ and $w_2$ are given by
\begin{eqnarray}
\label{eq_case1U_1}
U_1 = U_2  \sim  \frac{4\sqrt{2}}{\sigma} \frac{Z_1}{Z}.
\end{eqnarray}
Therefore, we derive $Z$ and $Z_1$ below.
The zeta function of $Z$ is given by
\begin{eqnarray}
\zeta(z) &=& \frac{1}{2 \pi}\int_{-\infty}^{\infty} \int_{-\infty}^{\infty} dw_1dw_2 (w_1^2w_2^2)^z \exp\left(-\frac{1}{2}(w_1^2+w_2^2)\right)\\
& = & \frac{1}{2 \pi}\int_{-\infty}^{\infty} dw_1 w_1^{2z} \exp\left(-\frac{w_1^2}{2}\right) \int_{-\infty}^{\infty} dw_2 w_2^{2z} \exp\left(-\frac{w_2^2}{2}\right) \\
& = & \frac{2}{\pi}\int_{0}^{\infty} dw_1 w_1^{2z} \exp\left(-\frac{w_1^2}{2}\right) \int_{0}^{\infty} dw_2 w_2^{2z} \exp\left(-\frac{w_2^2}{2}\right) \\
& = & \frac{2}{\pi}\int_{0}^{\infty} \frac{dx_1}{\sqrt{2 x_1}} (2 x_1)^{z} \exp\left(-x_1\right)  \int_{0}^{\infty} \frac{dx_2}{\sqrt{2 x_2}} (2 x_2)^{z} \exp\left(-x_2\right) \\
& = & \frac{2^{2z}}{\pi}\int_{0}^{\infty} dx_1 x_1^{z-\frac{1}{2}} \exp\left(-x_1\right)  \int_{0}^{\infty} dx_2 x_2^{z-\frac{1}{2}} \exp\left(-x_2\right) \\
&=& \frac{2^{2z}}{\pi}\Gamma\left(z+\frac{1}{2}\right)^2.
\end{eqnarray}
By applying the inverse Mellin transform to this equation, we can obtain the function $V(t)$ as follows:
\begin{eqnarray}
V(t)&=&\frac{1}{2\pi i}\int_{c-i\infty}^{c+i\infty}\zeta(z)t^{-z-1}dz\\
&=&\frac{1}{2\pi i}\int_{c-i\infty}^{c+i\infty}\frac{2^{2z}}{\pi}\Gamma\left(z+\frac{1}{2}\right)^2 t^{-z-1}dz\\
&=&\frac{1}{2\pi i}\int_{c-i\infty}^{c+i\infty}
    \frac{1}{\left( z + \frac{1}{2}\right)^2} \left\{\frac{2^{2z}}{\pi}\Gamma\left(z+\frac{3}{2}\right)^2 t^{-z-1} \right\}dz\\
&=& \left.\frac{d}{dz} \left\{ \frac{1}{t \pi} \left( \frac{2^2}{t}\right)^z \Gamma\left(z+\frac{3}{2}\right)^2 \right\} \right|_{z=-\frac{1}{2}}\\
&=& \frac{1}{2 \pi} t^{-\frac{1}{2}} \left\{2\log 2 - 2 \gamma - \log t\right\},
\end{eqnarray}
where $\gamma$ is the Euler constant.
Therefore, $Z$ can be obtained as follows:
\begin{eqnarray}
Z &=& \int_{0}^{\infty} dt \exp(-nt) V(t)\\
&=& \frac{1}{2 \pi} \int_{0}^{\infty} dt \exp(-nt) t^{-\frac{1}{2}} \left\{2\log 2 - 2 \gamma - \log t\right\}\\
&=& \frac{1}{2 \pi} \int_{0}^{\infty} \frac{dt'}{n} \exp(-t') \left(\frac{t'}{n}\right)^{-\frac{1}{2}} 
    \left\{\log n + 2\log 2 - 2 \gamma - \log t'\right\}\\
&=& \frac{1}{2\sqrt{\pi}} \frac{\log n}{n^{\frac{1}{2}}} \left\{ 1 + \frac{4 \log 2 - \gamma}{\log n} \right\}.
\end{eqnarray}\par
Next, we derive $Z_1$.
The zeta function of $Z_1$ is given by
\begin{eqnarray}
\zeta_1(z) &=& \frac{1}{2 \pi}\int_{-\infty}^{\infty} \int_{-\infty}^{\infty} dw_1dw_2 (w_1^2w_2^2)^z |w_1| 
    \exp\left(-\frac{1}{2}(w_1^2+w_2^2)\right)\\
& = & \frac{1}{2 \pi}\int_{-\infty}^{\infty} dw_1 w_1^{2z} |w_1| \exp\left(-\frac{w_1^2}{2}\right) 
    \int_{-\infty}^{\infty} dw_2 w_2^{2z} \exp\left(-\frac{w_2^2}{2}\right) \\
& = & \frac{2}{\pi}\int_{0}^{\infty} dw_1 w_1^{2z+1} \exp\left(-\frac{w_1^2}{2}\right) 
    \int_{0}^{\infty} dw_2 w_2^{2z} \exp\left(-\frac{w_2^2}{2}\right) \\
& = & \frac{2}{\pi}\int_{0}^{\infty} \frac{dx_1}{\sqrt{2 x_1}} (2 x_1)^{z+\frac{1}{2}} \exp\left(-x_1\right)  
    \int_{0}^{\infty} \frac{dx_2}{\sqrt{2 x_2}} (2 x_2)^{z} \exp\left(-x_2\right) \\
& = & \frac{2^{2z+\frac{1}{2}}}{\pi}\int_{0}^{\infty} dx_1 x_1^{z} \exp\left(-x_1\right)  
    \int_{0}^{\infty} dx_2 x_2^{z-\frac{1}{2}} \exp\left(-x_2\right) \\
&=&\frac{2^{2z+\frac{1}{2}}}{\pi} \Gamma(z+1)\Gamma\left(z+\frac{1}{2}\right).
\end{eqnarray}
By applying the inverse Mellin transform to this equation, we can obtain the function $V_1(t)$ as follows:
\begin{eqnarray}
V_1(t)&=&\frac{1}{2\pi i}\int_{c-i\infty}^{c+i\infty}\zeta_1(z)t^{-z-1}dz\\
&=&\frac{1}{2\pi i}\int_{c-i\infty}^{c+i\infty}\frac{2^{2z+\frac{1}{2}}}{\pi} \Gamma(z+1)\Gamma\left(z+\frac{1}{2}\right) t^{-z-1}dz\\
&=&\frac{1}{2\pi i}\int_{c-i\infty}^{c+i\infty}
    \frac{1}{ z + \frac{1}{2}} \left\{\frac{2^{2z+\frac{1}{2}}}{\pi}\Gamma(z+1)\Gamma\left(z+\frac{3}{2}\right) t^{-z-1} \right\}dz\\
&=& \left. \left\{ \frac{2^{2z+\frac{1}{2}}}{\pi}\Gamma(z+1)\Gamma\left(z+\frac{3}{2}\right) t^{-z-1} \right\} \right|_{z=-\frac{1}{2}}\\
&=&\frac{1}{\sqrt{2 \pi}} t^{-\frac{1}{2}}.
\end{eqnarray}
Therefore, $Z_1$ can be obtained as follows:
\begin{eqnarray}
Z_1 &=& \int dt \exp(-nt) V_1(t)\\
&=& \frac{1}{\sqrt{2 \pi}} \int dt \exp(-nt) t^{-\frac{1}{2}}\\
&=& \frac{1}{\sqrt{2 \pi}} \int \frac{dt'}{n} \exp(-t') \left(\frac{t'}{n}\right)^{-\frac{1}{2}}\\
&=& \frac{1}{\sqrt{2}}\frac{1}{n^{\frac{1}{2}}}.
\end{eqnarray}

By substituting these results into Eq. (\ref{eq_case1U_1}), we obtain $U_1$ as follows:
\begin{eqnarray}
U_1 &=& \frac{4\sqrt{2}}{\sigma} \frac{1}{\sqrt{2}}\frac{1}{n^{\frac{1}{2}}} 
    \left\{ \frac{1}{2\sqrt{\pi}} \frac{\log n}{n^{\frac{1}{2}}} \left\{ 1 + \frac{4 \log 2 - \gamma}{\log n} \right\}\right\}^{-1}\\
&=& \frac{8 \sqrt{\pi}}{\sigma (\log n + 4 \log 2 - \gamma)}.
\end{eqnarray}

\section{Analytical solution for $f(w_1,w_2) = w_1^2w_2^4$}
From Lemma 2, the average acceptance rate $U_1$ for $w_1$ is given by
\begin{eqnarray}
\label{eq_case2U_1}
U_1 & \sim & \frac{4\sqrt{2}}{\sigma} \frac{Z_1}{Z}.
\end{eqnarray}
Therefore, $Z$ and $Z_1$ are derived below.
For $Z$, its zeta function is given by
\begin{eqnarray}
\zeta(z) &=& \int_{-\infty}^{\infty} dw_1 \int_{-\infty}^{\infty} dw_2 f(w_1, w_2)^z\varphi(w_1, w_2)\\
 &=& \int_{-\infty}^{\infty} dw_1 \int_{-\infty}^{\infty} dw_2 (w_1^2 w_2^4)^z\frac{1}{2\pi}\exp\left(-\frac{1}{2}(w_1^2+w_2^2)\right)\\
 &=& \frac{2}{\pi} \int_{0}^{\infty} dw_1 w_1^{2z} \exp\left(-\frac{w_1^2}{2}\right) 
 \int_{0}^{\infty} dw_2 w_2^{4z} \exp\left(-\frac{w_2^2}{2}\right)\\ 
 &=& \frac{2}{\pi} \int_{0}^{\infty} \frac{dx_1}{\sqrt{2x_1}} (2x_1)^{z} \exp\left(-x_1\right)
 \int_{0}^{\infty} \frac{dx_2}{\sqrt{2x_2}} (2x_2)^{2z} \exp\left(-x_2\right)\\
 &=& \frac{2^{3z}}{\pi} \Gamma\left(z+\frac{1}{2}\right) \Gamma\left(2z + \frac{1}{2}\right).
\end{eqnarray}
By applying the inverse Mellin transform to this equation, 
we can obtain the function $V(t)$ as follows:
\begin{eqnarray}
V(t) &=& \frac{1}{2\pi i}\int_{c-i\infty}^{c+i\infty} \zeta(z) t^{-z-1}dz\\
 &=& \frac{1}{2\pi i} \int_{c-i\infty}^{c+i\infty} dz
    \frac{2^{3z}}{\pi} \Gamma\left(z+\frac{1}{2}\right) \Gamma\left(2z + \frac{1}{2}\right) t^{-z-1}\\
 &=& \frac{1}{2\pi i} \int_{c-i\infty}^{c+i\infty} dz
    \frac{1}{z+\frac{1}{4}} \left\{ \frac{2^{3z-1}}{\pi} \Gamma\left(z+\frac{1}{2}\right) \Gamma\left(2z + \frac{3}{2}\right) t^{-z-1} \right\}\\
 &=& \frac{2^{-\frac{7}{4}}}{\pi} \Gamma\left(\frac{1}{4}\right) t^{-\frac{3}{4}}.
\end{eqnarray}
Therefore, $Z$ can be obtained as follows:
\begin{eqnarray}
Z &=& \int_{0}^{\infty} dt \exp(-nt)V(t),\\
 &=& \frac{2^{-\frac{7}{4}}}{\pi} \Gamma\left(\frac{1}{4}\right)
 \int_{0}^{\infty} dt \exp(-nt) t^{-\frac{3}{4}},\\
 &=& \frac{1}{n^{\frac{1}{4}}} \frac{2^{-\frac{7}{4}}}{\pi} \left\{\Gamma\left(\frac{1}{4}\right)\right\}^2.
\end{eqnarray}
Next, we derive $Z_1$.
The zeta function of $Z_1$ is given by
\begin{eqnarray}
\zeta_1(z) &=& \int_{-\infty}^{\infty} dw_1 \int_{-\infty}^{\infty}  dw_2
    |w_1|f(w_1, w_2)^z \varphi(w_1, w_2)\\
 &=& \int_{-\infty}^{\infty} dw_1 \int_{-\infty}^{\infty} dw_2 
 |w_1|(w_1^2 w_2^4)^z\frac{1}{2\pi}
 \exp\left(-\frac{1}{2}(w_1^2+w_2^2)\right)\\
 &=& \frac{2}{\pi} \int_{0}^{\infty} dw_1 w_1^{2z+1} \exp\left(-\frac{w_1^2}{2}\right) 
 \int_{0}^{\infty} dw_2 w_2^{4z} \exp\left(-\frac{w_2^2}{2}\right)\\ 
 &=& \frac{2}{\pi} \int_{0}^{\infty} \frac{dx_1}{\sqrt{2x_1}} (2x_1)^{z+\frac{1}{2}} \exp\left(-x_1\right)
 \int_{0}^{\infty} \frac{dx_2}{\sqrt{2x_2}} (2x_2)^{2z} \exp\left(-x_2\right)\\
 &=& \frac{2^{3z+\frac{1}{2}}}{\pi} 
 \Gamma\left(z + 1\right) \Gamma\left(2z + \frac{1}{2}\right),
\end{eqnarray}
By applying the inverse Mellin transform to this equation, 
we can obtain the function $V_1(t)$ as follows:
\begin{eqnarray}
V_1(t) &=& \frac{1}{2\pi i}\int_{c-i\infty}^{c+i\infty} \zeta_1(z) t^{-z-1}dz,\\
 &=& \frac{1}{2\pi i}\int_{c-i\infty}^{c+i\infty} dz \frac{2^{3z+\frac{1}{2}}}{\pi} 
 \Gamma\left(z + 1\right) \Gamma\left(2z + \frac{1}{2}\right) t^{-z-1},\\  &=& \frac{1}{2\pi i}\int_{c-i\infty}^{c+i\infty} dz
 \frac{1}{z+\frac{1}{4}} \left\{\frac{2^{3z-\frac{1}{2}}}{\pi} 
 \Gamma\left(z + 1\right) \Gamma\left(2z + \frac{3}{2}\right) t^{-z-1}\right\},\\
 &=& \frac{2^{-\frac{5}{4}}}{\pi} \Gamma\left(\frac{3}{4}\right) t^{-\frac{3}{4}}.
\end{eqnarray}
Therefore, $Z_1$ can be obtained as follows:
\begin{eqnarray}
Z_1 &=& \int_{0}^{\infty} dt \exp(-nt)V_1(t),\\
 &=& \frac{2^{-\frac{5}{4}}}{\pi} \Gamma\left(\frac{3}{4}\right) 
    \int_{0}^{\infty} dt \exp(-nt) t^{-\frac{3}{4}},\\
 &=& \frac{1}{n^{\frac{1}{4}}} \frac{2^{-\frac{5}{4}}}{\pi} \Gamma\left(\frac{3}{4}\right)\Gamma\left(\frac{1}{4}\right).
\end{eqnarray}
By substituting these results into Eq. (\ref{eq_case2U_1}), we obtain $U_1$ as follows:
\begin{eqnarray}
 U_1 &=& \frac{Z_1}{Z}\times \frac{4\sqrt{2}}{\sigma},\\
 &=& \frac{8}{\sigma}\frac{\Gamma\left(\frac{3}{4}\right)}{\Gamma\left(\frac{1}{4}\right)}.
\end{eqnarray}

In the same way, we derive the average acceptance rate $U_2$ for $w_2$.
\begin{eqnarray}
\label{eq_case2U_2}
U_2 & \sim & \frac{4\sqrt{2}}{\sigma} \frac{Z_2}{Z}
\end{eqnarray}
Since $Z$ was obtained by deriving $U_1$, we derive only $Z_2$.
The zeta function of $Z_2$ can be obtained as follows:
\begin{eqnarray}
\zeta_2(z) &=& \int_{-\infty}^{\infty} dw_1 
    \int_{-\infty}^{\infty}  dw_2 
    |w_2| f(w_1, w_2)^z \phi(w_1, w_2)\\
 &=& \int_{-\infty}^{\infty} dw_1 \int_{-\infty}^{\infty} dw_2 
    |w_2|(w_1^2 w_2^4)^z\frac{1}{2\pi}
    \exp\left(-\frac{1}{2}(w_1^2+w_2^2)\right)\\
 &=& \frac{2}{\pi} \int_{0}^{\infty} dw_1 
    w_1^{2z} \exp\left(-\frac{w_1^2}{2}\right) 
    \int_{0}^{\infty} dw_2 w_2^{4z+1} \exp\left(-\frac{w_2^2}{2}\right)\\ 
 &=& \frac{2}{\pi} \int_{0}^{\infty} \frac{dx_1}{\sqrt{2x_1}} 
    (2x_1)^{z} \exp\left(-x_1\right)
    \int_{0}^{\infty} \frac{dx_2}{\sqrt{2x_2}} 
    (2x_2)^{2z+\frac{1}{2}} \exp\left(-x_2\right)\\
 &=& \frac{2^{3z+\frac{1}{2}}}{\pi} \Gamma\left(z+\frac{1}{2}\right) \Gamma\left(2z + 1\right).
\end{eqnarray}
By applying the inverse Mellin transform to this equation, 
we can obtain the function $V_2(t)$ as follows:
\begin{eqnarray}
V_2(t) &=& \frac{1}{2\pi i}\int_{c-i\infty}^{c+i\infty} \zeta_2(z) t^{-z-1}dz\\
 &=& \frac{1}{2\pi i}\int_{c-i\infty}^{c+i\infty} dz 
    \frac{2^{3z+\frac{1}{2}}}{\pi} \Gamma\left(z+\frac{1}{2}\right)
    \Gamma\left(2z + 1\right) t^{-z-1}\\
 &=& \frac{1}{2\pi i}\int_{c-i\infty}^{c+i\infty} dz 
    \frac{1}{\left(z+\frac{1}{2}\right)^2}
    \frac{2^{3z-\frac{1}{2}}}{\pi} \Gamma\left(z+\frac{3}{2}\right) \Gamma\left(2z + 2\right) t^{-z-1}\\
 &=& \left.\frac{d}{dz} \left\{ \frac{2^{3z-\frac{1}{2}}}{\pi}
    \Gamma\left(z+\frac{3}{2}\right) 
    \Gamma\left(2z + 2\right) t^{-z-1}\right\} \right|_{z=-\frac{1}{2}}\\
 &=& \frac{1}{4\pi} \left(3\log 2- 3 \gamma - \log t\right) t^{-\frac{1}{2}}.
\end{eqnarray}
Therefore, $Z_2$ can be obtained as follows:
\begin{eqnarray}
Z_2 &=& \int_{0}^{\infty} dt \exp(-nt)V_2(t),\\
 &=& \frac{1}{4\pi}\int_{0}^{\infty} dt \exp(-nt) \left(3\log 2- 3 \gamma - \log t\right) t^{-\frac{1}{2}},\\
 &=& \frac{\log n}{n^{\frac{1}{2}}} \frac{1}{4\pi} \Gamma\left(\frac{1}{2}\right)\left(1+\frac{\log 2 - 4 \gamma}{\log n}\right).
\end{eqnarray}
By substituting these results into Eq. (\ref{eq_case2U_2}), we obtain $U_2$ as follows:
\begin{eqnarray}
 U_2 &=& \frac{Z_2}{Z}\times \frac{4\sqrt{2}}{\sigma},\\
 &=& \frac{2^{\frac{9}{4}}}{\sigma}\frac{\log n}{n^{\frac{1}{4}}}\frac{\Gamma\left(\frac{1}{2}\right)}{\Gamma\left(\frac{1}{4}\right)^2}\left(1+\frac{\log 2 - 4\gamma}{\log n}\right).
\end{eqnarray}

\bibliographystyle{unsrt}  
\bibliography{references}

\end{document}